\documentclass[12pt]{article}
\usepackage{amsmath}
\usepackage{amssymb}
\usepackage{url}

\usepackage[T1]{fontenc}
\usepackage{mathptmx}
\usepackage{diagrams}
\diagramstyle{midshaft,PostScript=dvips,nohug}
\newarrow{Dotsto}.. ..>

\newtheorem{lemma}{Lemma}[section]
\newtheorem{prop}[lemma]{Proposition}
\newtheorem{thm}[lemma]{Theorem}

\newcommand{\A}{\mathcal{A}}
\newcommand{\E}{\mathcal{E}}

\newcommand{\oo}{\overline}
\newcommand{\Aff}{\operatorname{Aff}}

\newcommand{\Lin}{\operatorname{Lin}}
\newcommand{\spec}{\operatorname{spec}}

\title{Affine combinations in affine schemes}
\author{Anders Kock}
\date{}

\begin{document}
\maketitle

\section*{Introduction} The notion of ``neighbour points'' in 
algebraic geometry is a geometric rendering of the notion of 
nilpotent elements in commutative rings, and was deve\-loped since the 
time of Study, K\"{a}hler, Hjelmslev,  and notably in French algebraic 
geometry (Grothendieck, Weil et al.) since the 1950s. They introduced 
it via what they call {\em the first neighbourhood of the diagonal}.

In \cite{SDG}, \cite{DFIC} and \cite{SGM} the neighbour notion was considered on an 
axiomatic basis, essentially for finite dimensional manifolds; one of 
the aims was to describe a combinatorial theory of differential forms.

In the specific context of algebraic geometry, such theory of 
differential forms was developed 
in \cite{BM}, where it applies not only to manifolds, but to arbitrary 
schemes.

One aspect, present in \cite{DFIC} and \cite{SGM}, but not in \cite{BM}, is the possibility of 
forming affine combinations of finite sets of mutual (1st order) 
neighbour points. The present note completes this aspect, by giving 
the construction of such affine combinations, at least in the category 
of {\em affine} schemes (the dual of the category of finitely presented 
commutative rings or $k$-algebras).

The interest in having the possibility of such affine combinations is 
documented in several places in \cite{SGM}, and is in \cite{DFIC} the basis for 
constructing, for any manifold, a simplicial object, whose cochain 
complex is the deRham complex of the manifold.

From a more philosophical viewpoint, one may say that the possibility 
of having affine combinations, for sets of mutual neighbour points, 
expresses 
in a concrete way the idea that spaces are ``infinitesimally like affine 
spaces''.

\section{Neighbour  maps between algebras}\label{NMBAx}
Let $k$ be a commutative ring. Consider commutative $k$-algebras
$B$ and $C$ and two $k$-algebra maps $f$ and $g:B \to C$.
We say that they are {\em neighbours}, or more completely, 
{\em (first order) infinitesimal  neighbours}, if 
 \begin{equation}\label{nb2x} (f(a) - g(a))
\cdot (f(b)-g(b)) =0 \mbox{ for all } a,b \in B,\end{equation}
or equivalently, if\begin{equation}\label{nb1x} f(a)\cdot 
g(b)+g(a)\cdot f(b) = 
f(a\cdot b)+g(a\cdot b) \mbox{ for all } a,b \in B.\end{equation} 
(Note that this latter formulation makes no use of ``minus''.)
When this holds, we write $f\sim g$ (or more completely, $f\sim _{1}g$). The relation $\sim$ is a reflexive and 
symmetric relation (but not transitive).  If the element $2\in k$ is invertible, a third equivalent 
formulation of $f\sim g$ goes \begin{equation}\label{nb3x} 
(f(a)-g(a))^{2}=0 \mbox{ for all } a\in B.\end{equation}
For, it is clear that (\ref{nb2x}) implies (\ref{nb3x}). Conversely, 
assume (\ref{nb3x}), and let $a, b \in B$ be arbitrary, and apply 
(\ref{nb3x}) to the element $a+b$. Then by assumption, 
and using that $f$ and $g$ are 
algebra maps\footnote{``algebra'' means throughout ``commutative $k$-algebra'', and 
similarly for  algebra maps. When we say 
``linear map'', we mean $k$-linear.  By 
$\otimes$, we  mean $\otimes _{k}$.},
$$0=(f(a+b) -g(a+b))^{2} = [(f(a)-g(a))+(f(b)-g(b))]^{2}$$
$$= (f(a)-g(a))^{2}+(f(b)-g(b))^{2} - 2(f(a)-g(a))\cdot (f(b)-g(b)).$$
The two first terms are 0 by assumption, hence so is the third. Now 
divide by $2$.

\medskip 

Note that if $C$ has no zero-divisors, then $f\sim g$ is 
equivalent to $f=g$.

\medskip

It is clear that the relation $\sim$ is stable under precomposition:
\begin{equation}\label{precx}\mbox{ if }h:B'\to B \mbox{ and } 
f\sim g:B\to C\mbox{, then  }f\circ h \sim g\circ h: B' \to C,\end{equation}
and,  using that $h$ is an algebra map, it is also stable under postcomposition:
\begin{equation}\label{postcx}\mbox{ if }  h:C\to C'\mbox{ and }f\sim g:B\to C\mbox{, then  }h\circ f 
\sim h\circ g :B \to C'.\end{equation}
Also, if $h:B'\to B$ is a surjective algebra map, precomposition by 
$h$ not only preserves the neighbour relation, it also reflects it, 
in the following sense
\begin{equation}\label{reflxx}f\circ h \sim g\circ h \mbox{ implies } 
f\sim g.
\end{equation}
This is immediate from (\ref{nb2x}); the $a$ and $b$ occurring there 
is of the form $h(a')$ and $h(b')$ for suitable $a'$ and $b'$ in 
$B'$, by surjectivity of $h$.

  An alternative ``element-free'' formulation of the neighbour relation 
(Proposition \ref{annx} below) 
comes from a standard piece of commutative algebra. Recall that for 
commutative $k$-algebras $A$ and $B$, the tensor product $A\otimes B$ carries 
structure of commutative $k$-algebra ($A\otimes B$ is in fact a 
coproduct of $A$ and $B$); the multiplication map $m: B\otimes B \to B$ is a 
$k$-algebra homomorphism; so the kernel is an ideal $J\subseteq 
B\otimes B$.

The following is a classical description of  the ideal $J\subseteq 
B\otimes B$; we include it for completeness. 
\begin{prop}\label{polyx} The kernel $J$ of $m:B\otimes B \to B$ is 
  generated 
by the expressions $  1\otimes b - b\otimes 1$, for $b\in B$. Hence the 
ideal $J^{2}$ is generated by the expressions 
$(1\otimes a -a \otimes 1)\cdot (1\otimes b - b\otimes 1)$ (or 
equivalently, by the expressions 
$1\otimes 
ab  +ab\otimes 1-a\otimes b -b\otimes a$).
\end{prop}
{\bf Proof.} It is clear that $1\otimes b -b\otimes 1$ is in $J$. 
Conversely, assume that $\sum_{i}a_{i}\otimes b_{i}$ is in $J$, i.e.\ 
that $\sum_{i}a_{i}\cdot b_{i}=0$. Rewrite the $i$th term 
$a_{i}\otimes b_{i}$ as follows:
$$a_{i}\otimes b_{i}= a_{i}b_{i}\otimes 1 + (a_{i}\otimes 1)\cdot 
(1\otimes b_{i} - b_{i}\otimes 1)$$
and sum over $i$; since $\sum_{i}a_{i}b_{i}=0$, we are left with 
$\sum_{i} (a_{i}\otimes 1)\cdot 
(1\otimes b_{i} - b_{i}\otimes 1)$, which belongs to the $B\otimes B$-module 
generated by elements of the form  $1\otimes b - b\otimes 1$. -- The second assertion 
follows, since $ab\otimes 1+1\otimes 
ab -a\otimes b -b\otimes a$ is the product of the two 
generators $ 1\otimes a- a \otimes 1$ and $1\otimes b - b\otimes 1$.
(Note that the proof gave a slightly stronger result, namely that $J$ 
is generated already as a $B$-module, by the elements $  1\otimes b - 
b\otimes 1$, via the algebra map $i_{0}: B \to B\otimes B$, where  
$i_{0}(b)=b\otimes 1$).

\medskip

From the second assertion in this Proposition immediately follows 
that $f\sim g$ iff $\{f,g\}: B\otimes B \to C$ factors across the 
quotient map 
$B\otimes B \to (B\otimes B)/J^{2}$ (where $\{f,g\}: B\otimes B \to C$ 
denotes the map given by $a\otimes b \mapsto f(a)\cdot g(b)$); 
equivalently:
\begin{prop}\label{annx}For $f, g : B \to C$, we have $f\sim g$ if 
and only if $\{f,g\}:B\otimes B \to C$ annihilates $J^{2}$.
\end{prop}

\medskip 
  
 The two natural 
inclusion maps $i_{0}$ and $i_{1}: 
B \to B\otimes B$ (given by $b\mapsto b\otimes 1$ and $b\mapsto 
1\otimes b$, respectively) are not in general neighbours, but when 
postcomposed with $\pi: B\otimes B \to (B\otimes B )/J^{2}$, they are:
$$\pi\circ i_{0}\sim \pi\circ i_{1},$$
and this is in fact the universal pair of neighbour  algebra maps with domain 
$B$.

\section{Neighbours for polynomial algebras}\label{NPAx}
We consider the polynomial algebra $B:= k[X_{1}, \ldots ,X_{n}]$.
Identifying $B\otimes B$ with $k[Y_{1},\ldots,Y_{n},Z_{1}, \ldots 
,Z_{n}]$, the multiplication map $m$ is the algebra map given by $Y_{i}\mapsto 
X_{i}$ and $Z_{i}\mapsto X_{i}$, so it is clear that the kernel $J$ 
of $m$ contains the $n$ elements $Z_{i}-Y_{i}$. The following 
Proposition should be classical:
\begin{prop}\label{poly1x}The ideal $J \subseteq B\otimes B$, for $B=k[X_{1},\ldots 
,X_{n}]$, is generated (as a $B\otimes B$-module) by the $n$ elements 
$Z_{i}-Y_{i}$.
\end{prop}
{\bf Proof.} From Proposition \ref{polyx}, we know that $J$ is 
generated  by 
elements $P(\underline{Z})-P(\underline{Y})$, for $P\in 
k[\underline{X}]$ (where $\underline{X}$ denotes $X_{1}, \ldots 
,X_{n}$, and similarly for $\underline{Y}$ and $\underline{Z}$). So 
it suffices to prove that $P(\underline{Z})-P(\underline{Y})$ is of 
the form
$$\sum_{i=1}^{n}(Z_{i}-Y_{i})Q_{i}(\underline{Y},\underline{Z}).$$This 
is done by induction in $n$. For $n=1$, it suffices, by linearity, to prove this 
fact for each monomial $X^{s}$. And this follows from the identity
 \begin{equation}\label{indx}Z^{s}-Y^{s}=(Z-Y
)\cdot (Z^{s-1}+ 
Z^{s-2}Y + \ldots + ZY^{s-2}+ Y^{s-1})\end{equation}
(for $s\geq 1$; for $s=0$, we get $0$).
For the induction step:
Write $P(\underline{X})$ as a sum of increasing powers of $X_{1}$,
$$P(X_{1},X_{2}, \ldots )= P_{0}(X_{2}, \ldots )+X_{1}P_{1}( X_{2}, 
\ldots ) +X_{1}^{2}P_{2}(X_{2}, \ldots ).$$
Apply the induction hypothesis to the first term. The remaining 
terms are of the form $X_{1}^{s}Q_{s}(X_{2}, \ldots )$ with $s\geq 1$; then the 
difference to be considered is
$$Y_{1}^{s}Q_{s}(Y_{2},\ldots) - Z_{1}^{s}Q_{s}(Z_{2}, \ldots )$$
which we may write as
$$Y_{1}^{s}(Q_{s}(Y_{2},\ldots)-Q_{s}(Z_{2},\ldots))+ 
Q_{s}(Z_{2},\ldots)(Y_{1}^{s}-Z_{1}^{s}).$$
The first term in this sum is taken care of by the induction 
hypothesis, the second term uses the identity (\ref{indx})
which shows that this term is in the ideal generated by $(Z_{1}-Y_{1})$.

\medskip

From this follows immediately
\begin{prop}\label{poly2x}The ideal $J^{2}\subseteq B\otimes B$, for $B= k[X_{1}, \ldots ,X_{n}]$ is 
generated (as a $B\otimes B$-module) by the elements $(Z_{i}-Y_{i})(Z_{j}-Y_{j})$ (for $i,j= 1, 
\ldots ,n$) (identifying $B\otimes B$ with $k[Y_{1}, \ldots ,Y_{n}, 
Z_{1}, \ldots ,Z_{n}]$).
\end{prop} 
\medskip

(The algebra $(B\otimes B)/J^{2}$ is the algebra representing the 
affine scheme 
``first neighbourhood of the diagonal'' for the affine scheme 
represented by $B$, alluded to in the introduction.)

\medskip

Algebra maps $\underline{a}:k[X_{1}, \ldots ,X_{n}] \to C$ are 
completely given by an $n$-tuple of ele\-ments $a_{i}:=  
\underline{a}(X_{i})\in C$ ($i=1, \ldots ,n$). Let $\underline{b}: k[X_{1}, 
\ldots ,X_{n}] \to C$ be similarly given by the $n$-tuple $b_{i}\in 
C$. The  decision when $\underline{a}\sim \underline{b}$ can be 
expressed equationally in terms of these two $n$-tuples of elements 
in $C$, i.e.\ as a purely equationally described condition on 
elements $(a_{1}, 
\ldots ,a_{n}, b_{1}, \ldots ,b_{n}) \in C^{2n}$:

\begin{prop}\label{Dnx} Consider two algebra maps $\underline{a}$ and 
$\underline{b}: k[X_{1}, \ldots ,X_{n}] \to C$. Let 
$a_{i}:=\underline{a}(X_{i})$ and $b_{i}:=\underline{b}(X_{i} )$. 
Then we have $\underline{a}\sim \underline{b}$ if and only if
\begin{equation}\label{nbbx}(b_{i}-a_{i})\cdot 
(b_{j}-a_{j})=0\end{equation}
for all $i,j = 1, \ldots ,n$. 
\end{prop}
{\bf Proof.} We have that 
$\underline{a}\sim \underline{b}$ iff the algebra map 
$\{\underline{a},\underline{b} \}$ annihilates the ideal $J^{2}$ 
for the algebra $k[X_{1}, \ldots ,X_{n}]$; and this in turn is equivalent to that it 
annihilates the set of generators for $J^{2}$ described in the 
Proposition \ref{poly2x}). But 
$\{\underline{a},\underline{b}\}((Z_{i}-Y_{i})\cdot (Z_{j}-Y_{j})) = 
(b_{i}-a_{i})\cdot (b_{j}-a_{j})$, and 
then the result is immediate.

\medskip

We therefore also say that the pair of $n$-tuples of elements in $C$ 
$$\left[ \begin{array}{ccc}a_{1}& 
\ldots &a_{n}\\  b_{1}& \ldots &b_{n}\end{array}\right]$$ are {\em 
neighbours} if 
(\ref{nbbx}) holds. 

For brevity, we call an $n$-tuple $(c_{1}, \ldots ,c_{n})$
of elements in $C^{n}$ a {\em 
vector}, and denote it $\underline{c}$ . Thus a vector $(c_{1}, 
\ldots ,c_{n})$ is neighbour of the 
``zero'' vector $\underline{0}=(0,\ldots,0)$ iff $c_{i}\cdot c_{j}=0$ for all $i$ 
and $j$.

\medskip

\noindent{\bf Remark.} Even when $2\in k$ is invertible, one cannot 
conclude that $\underline{a}\sim \underline{b}$ follows from
$(b_{i}-a_{i})^{2} =0$ for all $i=1, \ldots ,n$. 
For, consider $C:= k[\epsilon_{1}, \epsilon_{2}] = k[\epsilon] 
\otimes k[\epsilon]$ (where $k[\epsilon ]$ is the ``ring of dual 
numbers over $k$'', so $\epsilon^{2}=0$). Then the pair of $n$-tuples 
($n=2$ here) given by $(a_{1}, a_{2})= (\epsilon_{1}, \epsilon 
_{2})$ and $(b_{1},b_{2}) :=(0,0)$
has $(a_{i}-b_{i})^{2} = \epsilon_{i}^{2}=0$ for $i=1,2$, but 
$(a_{1}-b_{1})\cdot (a_{2}-b_{2}) = \epsilon_{1}\cdot \epsilon_{2} $, 
which is not $0$ in $C$.

\medskip

We already have the notion of when  two algebra maps $f$ and $g:B\to C$ are 
neighbours, or infinitesimal neighbours.  We also say that the pair 
$(f,g)$ form an {\em infinitesimal 1-simplex} (with $f$ and $g$ as 
{\em vertices}). Also, we have the 
derived (\ref{nbbx}) notion of when two vectors in $C^{n}$ 
  are 
neighbours, or form an infinitesimal 1-simplex. This terminology is 
suited for being generalized to defining the notion of {\em infinitesimal 
$p$-simplex} of algebra maps $B\to C$, or of {\em infinitesimal 
$p$-simplex} of vectors in $C^{n}$ (for $p=1,2, \ldots$).

Proposition \ref{Dnx}  generalizes immediately to infini\-te\-simal 
$p$-simpli\-ces (where the Proposition is the special case of $p=1$):

\begin{prop}\label{Dnpx} Consider $p+1$ algebra maps 
$\underline{a}_{i}:k[X_{1},\ldots ,X_{n}]\to C$ (for $i=0, \ldots ,p$), 
and let $a_{ij}\in C$ be 
$\underline{a}_{i} (X_{j})$, for $j=1, \ldots n$.
 Then the $\underline{a}_{i}$ form an 
infinitesimal $p$-simplex iff for all $i, i'=0, \ldots p$ and  $ 
j,j'=1, \ldots ,n$
\begin{equation}\label{Dtildex}(a_{ij}-a_{i'j})\cdot (a_{ij'}-a_{i'j'})=0.
\end{equation} 
\end{prop}

\section{Affine combinations of mutual neighbours}\label{ACMNx}

Let $C$ be a $k$-algebra. An {\em affine} combination in a $C$-module means 
here a linear combination in the module,  with coefficients from $C$, 
and where the 
sum of the coefficients is $1$ . 
We consider in particular the $C$-module $\Lin_{k}(B,C)$ of $k$-linear maps $B\to C$, 
where $B$ is another  $k$-algebra.
Linear combinations of algebra maps are 
linear, but may fail to preserve the multiplicative structure and $1$.
 However
\begin{thm}\label{onex}Let $f_{0}, \ldots ,f_{p}$ be a $p+1$-tuple of mutual 
neighbour algebra maps $B \to C$, and let $t_{0}, \ldots ,t_{p}$ be elements 
of $C$ 
with $t_{0}+\ldots +t_{p}=1$. Then the affine combination
$$\sum_{i=0}^{p}t_{i}\cdot f_{i}:B \to C$$
is an algebra map. Composing with a map $h: C\to C'$ preserves the 
affine combination.
\end{thm}
{\bf Proof.} Since the sum is a $k$-linear map, it suffices to prove that 
it preserves the multiplicative structure. It clearly preserves 1. To 
prove that it preserves products $a\cdot b$, we should compare
$$(\sum_{i}t_{i}f_{i}(a))\cdot(\sum_{j}t_{j}f_{j}(b)) =\sum 
_{i,j}t_{i}t_{j}f_{i}(a)\cdot f_{j}(b)$$
with
$\sum t_{i}f_{i}(a\cdot b)$.
Now use that $\sum_{j}t_{j}=1$; then $\sum t_{i}f_{i}(a\cdot b)$ may be rewritten as
$$\sum_{ij}t_{i}t_{j}f_{i}(a\cdot b).$$
Compare the two displayed double sums: the terms with 
$i=j$ match since each $f_{i}$ preserves multiplication. Consider a pair 
of indices $i\neq j$; the  terms with index $ij$ and $ji$ from the 
first sum contribute $t_{i}t_{j}$ times 
\begin{equation}\label{sux1}f_{i}(a)\cdot f_{j}(b)  + f_{j}(a)\cdot 
f_{i}(b),\end{equation}
and the terms terms with index $ij$ and $ji$ from the second sum 
contribute $t_{i}t_{j}$ times 
\begin{equation}\label{sux2}f_{i}(a\cdot b) + f_{j}(a\cdot 
b),\end{equation}
and the two displayed contributions are equal, since $f_{i}\sim 
f_{j}$ (use the formulation (\ref{nb1x})).
The last assertion is obvious from the construction.
\medskip
\begin{thm}\label{twox}Let Let $f_{0}, \ldots ,f_{p}$ be a $p+1$-tuple of mutual 
neighbour algebra maps $B \to C$. Then any two affine combinations (with 
coefficients from $C$)  of these 
maps are neighbours. 
\end{thm}
{\bf Proof.} Let $\sum _{i}t_{i}f_{i}$ and $\sum _{j}s_{j}f_{j}$ be 
two such affine combinations. To prove that they are neighbours means 
(using (\ref{nb1x})) 
to prove that for all $a$ and $b$ in $B$,  
$$(\sum_{i}t_{i}f_{i}(a))\cdot (\sum_{j}s_{j}f_{j}(b)) + 
(\sum_{j}s_{j}f_{j}(a))\cdot (\sum_{i}t_{i}f_{i}(b))$$
equals
$$\sum_{i}t_{i}f_{i}(a\cdot b) + \sum _{j}s_{j}f_{j}(a\cdot b).$$
The first of these expressions equals
$$\sum_{ij}t_{i}s_{j}f_{i}(a)\cdot f_{j}(b) + \sum _{ij}t_{i}s_{j}f_{j}(a)\cdot 
f_{i}(b)= \sum _{ij}t_{i}s_{j}[f_{i}(a)\cdot f_{j}(b)+ f_{j}(a)\cdot f_{i}(b)]$$
For the second expression, we use $\sum_{j}s_{j}=1$ and $\sum_{i}t_{i}=1$, to 
rewrite it as the left hand expression in
$$\sum_{ij}t_{i}s_{j}f_{i}(a\cdot b)  + 
\sum_{ij}t_{i}s_{j}f_{j}(a\cdot b)= \sum_{ij}t_{i}s_{j}[f_{i}(a\cdot 
b) + f_{j}(a\cdot b)].$$
For each $ij$, the two square bracket expression match by (\ref{nb1x}), since $f_{i}\sim 
f_{j}$.

\medskip

Combining these two results, we have

\begin{thm}\label{threex} Let $f_{0}, \ldots ,f_{p}$ be a $p+1$-tuple of mutual 
neighbour algebra maps $B \to C$. Then in the $C$-module of $C$-linear 
maps $B\to C$, the affine subspace $\Aff _{C}(f_{0},\ldots ,f_{p})$ 
 of affine combinations (with coefficients from $C$) of the 
$f_{i}$s  consists of algebra maps, and they are mutual neighbours.
\end{thm}


 In \cite{BM}, they describe  an ideal 
 $ J^{(2)}_{0p}$. It is the sum of  ideals 
$J_{rs}^{2}$ in the
$p+1$-fold tensor product $B\otimes \ldots \otimes B$, 
where $J_{rs}$ is the ideal generated by $i_{s}(b) - 
i_{r}(b)$ for $b\in B$ and $r<s$. We shall here denote it just 
$\overline{J}^{(2)}$ for brevity;   it has the property that the 
$p+1$ 
inclusions $B\to B\otimes \ldots \otimes B)$ become mutual 
neighbours, when composed with the quotient map $\pi: B\otimes \ldots 
\otimes B\to (B\otimes \ldots \otimes B)/\overline{J}^{(2)}$, and this is in 
fact the universal $p+1$ tuple of mutual neighbour maps with 
domain $B$.

We may, for any given $k$-algebra $B$,  encode  the construction of Theorem \ref{onex} into one single 
canonical map which does not mention any individual $B \to C$, by 
using the 
universal $p+1$-tuple, and the 
generic $p+1$ tuple of coefficients with sum 1, meaning $(X_{0}, 
X_{1}, \ldots ,X_{p})\in k[X_{1}, \ldots ,X_{p}]$ (where $X_{0}$ 
denotes $1-(X_{1}+ \ldots +X_{p})$;
namely as a $k$-algebra map
\begin{equation}\label{axxx}B\to (B^{\otimes k+1}/\overline{J}^{(2)} )\otimes k[X_{1}, 
\ldots ,X_{p}].\end{equation}
For, by the Yoneda Lemma, this is equivalent to giving
a (set theoretical) map, natural in $C$,
$$\hom((B^{\otimes p+1}/\overline{J}^{(2)} )\otimes k[X_{1}, \ldots 
,X_{p}], C) \to \hom (B,C),$$
(where $\hom$ denotes the set of $k$-algebra maps). An element on the 
left hand side is given by a $p+1$-tuple of mutual neighbouring 
algebra  maps $f_{i}: B\to C$, together with  
 a $p$-tuple $(t_{1}, \ldots ,t_{p})$ of 
elements in $C$. With $t_{0}:=1- \sum_{1}^{p} t_{i}$, such data 
produce an element $\sum_{0}^{p} t_{i}\cdot f_{i}$ in $\hom(B,C)$, by Theorem 
\ref{onex},and the construction is natural in $C$ by the last assertion in the 
Theorem.  

The affine scheme defined by the algebra $B^{\otimes 
p+1}/\overline{J}^{(2)}$ is (essentially) called $\Delta ^{(p)}_{B}$ in 
\cite{BM}, and, (in axiomatic context, and for manifolds, in a 
suitable sense), the corresponding object  is 
called $M_{[p]}$ in 
\cite{DFIC} and $M_{(1,1, \ldots,1)}$ in \cite{SDG} I.\ 18 (for suitable $M$). 

\section{Affine combinations in a $k$-algebra $C$}\label{ACkCx}
The constructions and results of the previous Section concerning
 infinitesimal $p$-simplices of 
algebra maps $B\to C$, specializes (by taking $B= k[X_{1}, \ldots 
,X_{n}]$, as in Section \ref{NPAx}) to infinitesimal $p$-simplices of vectors in
$C^{n}$; such a $p$-simplex is conveniently exhibited in a $(p+1)\times 
n$ matrix with entries $a_{ij}$ from $C$:
$$\left[ \begin{array}{ccc}
a_{01}& \ldots&a_{0n}\\
a_{11}& \ldots & a_{1n}\\
\vdots & & \vdots\\
a_{p1}&\ldots &a_{pn}
\end{array}\right].$$
We may of course form affine (or even linear) combinations, with 
coefficients from $C$, of the rows of this matrix, 
whether or not the rows are mutual neighbours.  
 But the same affine combination of 
the corresponding algebra maps  is in general only a $k$-linear 
map, not an algebra map. However, if the rows  are mutual neighbours
 in $C^{n}$, and hence the corresponding algebra maps are mutual 
neighbouring  algebra maps $k[X_{1}, \ldots ,X_{n}]\to C$, we have, by 
Theorem \ref{onex} that the affine combinations of the rows of the 
matrix corresponds to the similar affine combination of the algebra 
maps. For, it suffices to check their equality on the $X_{i}$s, since 
the $X_{i}$s generate $k[X_{1}, \ldots ,X_{n}]$ as an algebra. Therefore, the Theorems \ref{twox} and  \ref{threex} 
immediately translate into theorems about $p+1$-tuples of mutual 
neighbouring $n$-tuples of elements in the algebra $C$; recall that 
such a $p+1$-tuple may be identified with the rows of a $(p+1)\times 
n$ matrix with entries from $C$, satisfying the equations 
(\ref{Dtildex}).
\begin{thm}\label{TWOx} Let the rows of a $(p+1)\times n$ matrix with 
entries from $C$  be 
mutual neighbours. Then any two affine combinations (with 
coefficients from $C$) of these rows are neighbours. The set of all 
such affine combinations form an affine subspace of the $C$-module 
$C^{n}$. 
\end{thm}


Let us consider in particular the case where the $0$th row of a 
$(p+1)\times n$ matrix is the zero vector $(0, \ldots ,0)$. Then the following is an 
elementary calculation:

\begin{prop}Consider a $(p+1)\times n$ matrix $\{a_{ij}\}$ as above, but 
with $a_{0j}=0$ for $j=1, \ldots n$. Then the rows form an 
infinitesimal $p$-simplex iff 
 the conjunction of
\begin{equation}\label{ggx}a_{ij}\cdot a_{i'j'}+a_{i'j}\cdot 
a_{ij'}=0\mbox{ for all } i,i'=1, \ldots p, j,j'=1, \ldots n.\end{equation}
hold. 
and \begin{equation}\label{ffx}
a_{ij}\cdot a_{ij'}=0 \mbox{ for all }i=1, \ldots ,p, j= 1, \ldots 
n\end{equation}
 If $2$ is invertible in $C$, the equations (\ref{ffx}) follow from 
(\ref{ggx}).
\end{prop}
{\bf Proof.}
The last assertion follows by putting $i=i'$ in (\ref{ggx}), and 
dividing by $2$.
Assume that the rows of the matrix form an infinitesimal 
$p$-simplex. Then (\ref{ffx}) follows from 
$\underline{a}_{i}\sim \underline{0}$. The equation which asserts that 
$\underline{a}_{i}\sim \underline{a}_{i'}$ (for $i, i'=1, \ldots ,p$) 
is
$$(a_{ij}-a_{i'j})\cdot (a_{ij'}-a_{i'j'})=0\mbox{ for all } j,j'=1, 
\ldots n.$$
Multiplying out gives four terms, two of which vanish by virtue of 
(\ref{ffx}), and the two remaining add up to (minus) the sum on the 
left of (\ref{ggx}). For the converse implication, (\ref{ffx}) give that 
the last $p$ rows are $\sim \underline{0}$;  and (\ref{ffx}) and 
(\ref{ggx}) jointly give that $\underline{a}_{i}\sim 
\underline{a}_{i'}$, by essentially the same calculation which we 
have already made.

\medskip 

When $\underline{0}$ is one of the vectors in a $p+1$-tuple, any 
linear combination of the remaining $p$ vectors has the same value as 
a certain affine combination of all $p+1$ vectors, since the coefficient 
for $\underline{0}$ may be chosen arbitrarily without changing the 
value of the linear combination. Therefore the results on affine 
combinations of the rows in the $(p+1)\times n$ matrix with 
$\underline{0}$ as top row immediately translate to results about 
linear combinations of the remaining rows, i.e.\ they translate into 
results about $p\times n$ matrices, satisfying the equations 
(\ref{ggx}) and (\ref{ffx}); even the equations (\ref{ggx}) suffice, 
if $2$ is invertible.
In this form, the results were obtained in the preprint \cite{MNE}, and 
are stated here for completeness. We assume that $2\in k$ is 
invertible.  

We use the notation from 
\cite{SDG} I.16 and I. 18, where
 set of $p\times n$ matrices $\{a_{ij}\}$ satisfying 
(\ref{ggx}) was denoted $\tilde{D}(p,n)\subseteq C^{p\cdot n}$ (we there consider 
algebras $C$ over $k={\mathbb Q}$, so (\ref{ffx}) 
follows). In particular $\tilde{D}(2,2)$ consists of matrices of the 
form
$$\left[ \begin{array}{cc}
a_{11}&a_{12}\\
a_{21}&a_{22}
\end{array}\right] \mbox{ \quad with \quad }a_{11}\cdot 
a_{22}+a_{12}\cdot a_{21}=0.$$
Note that the determinant of such a matrix is $2$ times the product of 
the diagonal entries. And also note that $\tilde{D}(2,2)$ is stable 
under transposition of matrices.

The notation $\tilde{D}(p,n)$ may be consistently augmented to the case where 
$p=1$; we say $(a_{1}, \ldots ,a_{n}) \in \tilde{D}(1,n)$ if it is 
neighbour of $\underline{0}\in C^{n}$, i.e.\ if $a_{j}\cdot a_{j'} 
=0$ for all $j,j'=1, \ldots n$. (In \cite{SDG}, $\tilde{D}(1,n)$ is also 
denoted $D(n)$, and $D(1)$ is denoted $D$.)

It is clear that a $p\times n$ matrix belongs to $\tilde{D}(p,n)$ 
precisely when all its $2\times 2$ sub-matrices do; this is just a 
reflection of the fact that the defining equations (\ref{ggx}) only 
involve two row indices and two column indices at a time. From the 
transposition stability of $\tilde{D}(2,2)$ therefore follows that 
transposition $p\times n$ matrices takes $\tilde{D}(p,n)$ into 
$\tilde{D}(n,p)$.

Note that each of the rows of a matrix in $\tilde{D}(p,n)$ is a neighbour of 
$\underline{0}\in C^{n}$.

The results about affine combinations now get the following 
formulations in terms of linear combinations of the rows of matrices 
in $\tilde{D}(p,n)$:

\begin{thm}Given a matrix $X\in \tilde{D}(p,n)$. Let a $(p+1)\times 
n$  matrix $X'$ 
be obtained by adjoining to $X$ a row which is a linear combination of 
the rows of $X$. Then $X'$ is in $\tilde{D}(p+1,n)$. 
\end{thm}

\section{Geometric meaning} This Section contains essentially only reformulations 
of the previous material into geometric terms, and thereby it also 
contains  some 
motivation of the notions.

In any category $\E$, a map $f: I \to M$ may be thought of either as an 
$I$-para\-metrized family of elements (or points) of $M$, or as an 
$M$-valued function on $I$. This terminology is convenient in 
particular when $\E$ is in some sense a category of spaces, and is 
even traditional in algebraic geometry, for instance when $\E$ is the 
{\em dual} of the category $\A$ of commutative $k$-algebras. The 
category $\E$ is in this case essentially the category of {\em affine 
schemes} over $k$. We elaborate a little on this terminology for this 
specific case. When a commutative $k$-algebra $B$ is seen in the dual 
category $\E$, it is often  denoted $\spec{B}$ or  $\oo{B}$. If $X$ 
is an object in $\E$, the corresponding algebra is often denoted 
$O(X)$, and called the {\em algebra of functions} on $X$; more 
precisely, the algebra of {\em scalar} valued functions on $X$.

The category $\E$ is in this case a category equipped with a 
canonical commutative ring object $R$, namely $\oo{k[X]}$, whose 
geometric meaning is that it is the number line. The ring structure 
of $R$ in $\E$ comes about from the canonical co-ring structure of 
$k[X]$ in the category of commutative $k$-algebras. Now ring 
structure on the geometric line is elementary and well understood, since the time of 
Euclid, essentially, whereas the notion of coring is not  elementary, 
and is a much more recent invention. This is why $(\E ,R)$ is well 
suited to axiomatic abstraction, as in \cite{SDG}.

The reason why $\hom_{\E}(I,R)$ is a commutative ring (even a 
$k$-algebra) in $\E$ is that it is isomorphic to $O(I)$; for,
\begin{equation}\label{scax}\hom_{\E}(I,R) \cong \hom_{\E}(I, \oo{k[X]}) \cong 
\hom_{\A}(k[X],O(I)) \cong O(I),\end{equation}
the last isomorphism because $k[X]$ is the free $k$-algebra in one 
generator $X$.

Having a commutative ring object like $R$ in a category $\E$ is the first necessary 
condition for having the wonderful tool of coordinates available for 
the geometry in $\E$. 

Thus we have the {\em coordinate vector spaces} $R^{n}$; in the 
category of affine sche\-mes, this is $\oo{k[X_{1}, \ldots ,X_{n}}]$.
This object is an $R$-module object. The $R$-module structure may be 
described in the same way, in terms of $\oo{C}$-parametrized points of 
$R^{n}$, 
as when we described the ring structure of $R$ in such terms; 
a $\oo{C}$-parametrized point of $R^{n}$ amounts to an $n$-tuple of 
elements in the algebra $C$, and $C$-parametrized points of $R$ 
amount to elements of $C$. So the $R$-module structure of $R^{n}$ 
comes about from the $C$-module structure of $C^{n}$, for arbitrary 
$C$.

If $C$ is a finitely presented $k$-algebra, the object (space) 
$\oo{C}$ embeds into some $R^{n}$, since a finite presentation of $C$, 
with $n$ generators gives rise to a surjective (hence epimorphic) algebra map $k[X_{1}, 
\ldots ,X_{n}] \to C$.

Since we have the notion of ``neighbours'' for algebra maps, we have a 
notion of neighbours for maps in $\E$; and  it is preserved by pre- 
and post-composition. In particular, the embedding map $e: \oo{C} \to 
R^{n}$, 
obtained from a presentation of $C$ with $n$ generators, preserves the 
property of being neighbours, for parametrized families of 
points of $\oo{C}$. The embedding also {\em reflects} the neighbour 
relation, in the sense that if $x$ and $y$ are  
points of $\oo{C}$,  and $e(x)\sim e(y)$, 
then $x\sim y$. This is 
just a reformulation of (\ref{reflxx}).

Note the 
traditional replacement of the notation $e\circ x$ by $e(x)$. This is 
the ``symbolic'' counterpart of considering maps $I\to M$ as 
($I$-parametrized) {\em 
points}
\footnote{In algebraic geometry,  the terminology ``$I$-valued point 
of $M$'' is also used, see e.g.\ \cite{BM} p.\ 209. In \cite{SDG}, 
Part II, such a thing is called a ``generalized element of $M$, defined at 
stage $I$'', and a more elaborate description of the `logic' of 
generalized elements is presented.
}
 of $M$. Furthermore, it is tradition in algebraic geometry not 
always to be specific about the space $I$ of parameters for a 
parametrized point $I\to M$; thus, one talks about ``points $(x,y)$ of the 
unit circle $S$ given by $x^{2}+y^{2}=1$'', without explicit mention of whether it 
one means a real point, a rational point, a complex point, \ldots ; 
therefore, $I$ is omitted from notation, and one writes $(x,y)\in S$. 
This in particular applies when the statement or notion applies to 
{\em any} (parametrized) point of $M$, regardless of its parameter 
space $I$.

An example of this usage is for the neighbour relation in affine 
schemes. Consider  two  map $f$ and $g$ between $k$-algebras $B$ and $C$, as in 
(\ref{nb2x}). In the category of affine schemes, these are then 
neighbour maps $\oo{f}$ and $\oo{g}$:  $\oo{C}\to \oo{B}$,  i.e.\ 
neighbour points of $\oo{B}$ (parametrized by $\oo{C}$); we write 
$\oo{f}\in \oo{B}$, $\oo{g} \in \oo{B}$.

With this usage, Proposition \ref{Dnx} may be reformulated as
\begin{prop}Given two points $(a_{1}, \ldots ,a_{n})$ and $(b_{1}, 
\ldots ,b_{n}) \in R^{n}$. Then they are neighbours iff 
\begin{equation}\label{nbbbx}(b_{i}-a_{i})\cdot 
(b_{j}-a_{j})=0\end{equation}
for all $i,j = 1, \ldots ,n$. 
\end{prop}
Here the (common) parameter space  $\oo{C}$ of the $a_{i}s$ and $b_{i}$s is not 
mentioned explicitly; it could be any affine scheme. Note that 
(\ref{nbbbx}) is typographically the same as (\ref{nbbx}); in 
(\ref{nbbbx}), $a_{i}$ and $b_{j}$ are (parametrized) points of $R$ 
(parametrized by $\oo{C})$, in (\ref{nbbx}), they are elements in the 
algebra $C$; but these data correspond, by  (\ref{scax}), and this 
correspondence preserves algebraic structure. 

\medskip

Similarly, Proposition \ref{Dnpx} gets the reformulation:
\begin{prop}A $p+1$-tuple $\{ a_{ij}\}$ of points in $R^{n}$ form an infinitesimal 
$p$-simplex iff the equations (\ref{Dtildex}) hold. 
\end{prop}

This formulation, as the other formulations in ``synthetic'' terms, 
are the ones that are suited to axiomatic treatment, as in Synthetic 
Differential Geometry, which almost exclusively\footnote{Exceptions 
are found in \cite{L:D} (where $R$ is {\em constructed} out of an 
assumed infinitesimal object $T$); and in \cite{OPSC} and\cite{END}, where part of the 
reasoning does not assume any algebraic notions.} assumes a given commutative ring object $R$ in a category 
$\E$, preferably a topos, as a basic 
ingredient in the axiomatics. (The category $\E$ of affine schemes is not a 
topos, but the category of presheaves on $\E$ is, and it, and some of 
its subtoposes, are the basic categories considered in modern 
algebraic geometry, like \cite{DG}.)

\medskip

As a further illustration of the ``synthetic'' language, the 
algebraic formulation of the neighbour relation between algebra maps 
given in (\ref{nb3x}) (assuming $2\in k$ is invertible) may be 
rendered:

\begin{prop}For any affine scheme $\oo{B}$, scalar valued functions 
on $\oo{B}$ detect when points $x$ and $y$ of $\oo{B}$ are neighbours; 
i.e.\ if $\alpha (x) \sim \alpha (y)$ for all $\alpha : \oo{B}\to R$, 
then $x\sim y$. 
\end{prop}
Here $x$ and $y$ are points of $\oo{B}$, say parametrized by 
$\oo{C}$, i.e.\ they are maps $\oo{C}\to \oo{B}$, so they correspond to algebra 
maps $f$ and $g:B \to C$;  and $\alpha :\oo{B} \to R$ corresponds to 
$a\in B$.

\medskip

The fact  (\ref{reflxx}) gets the 
following formulation:
\begin{prop}\label{fivefx}Given  an affine scheme $\oo{B}$. Then for 
any finite presentation (with $n$ generators, say) of the algebra, 
the corresponding embedding $e: \oo{B} \to R^{n}$ (preserves and) 
reflects the relation $\sim$. 
  \end{prop}

So for $x$ and $y$ points in $R^{n}$, 
they come via $e$ from a pair of neighbour points in $\oo{B}$ iff 
they satisfy 1) the equations in $n$ variables defining $B$ in the 
presentation; and 2) satisfy the equations for being neighbours in 
$R^{n}$. In other words, the intrinsically defined neighbour relation on 
$\oo{B}$ (essentially (\ref{nb2x})) may be described purely equationally, using an finite equational 
presentation of $B$.

Or, in more elementary tems, which the synthetic tradition is very 
apt for utilizing:
Given a finite set of equations with coefficients from $R$. If 
of $x_{0},\ldots,x_{p}$ are points in $R^{n}$  and each of them 
satisfies the 
equations, then so does any affine combination of them, provided the 
points  
are mutual neighbours.
But note that Proposition \ref{fivefx}, together with the 
constructions of Sections \ref{NMBAx} and 
\ref{ACMNx} allow us to 
conclude that the neighbour conditions, and the point constructed by 
affine combinations in $R^{n}$, is intrinsic to the affine scheme $\oo{B}$ in question, and does not 
depend on an equational presentation of $B$.

\small

\noindent Anders Kock,
August 2015\\
Dept.\ of Mathematics,
University of Aarhus, Denmark\\
kock (at) math.au.dk

\end{document}